\documentclass[reqno]{amsart}
\usepackage{amssymb}
\usepackage{cite}
\allowdisplaybreaks
\date{\today}

\theoremstyle{plain}
\newtheorem{thm}{Theorem}[section]

\newtheorem{lem}[thm]{Lemma}
\newtheorem{prop}[thm]{Proposition}
\theoremstyle{definition}

\theoremstyle{remark}
\newtheorem{rem}{Remark}[section]
\numberwithin{equation}{section}

\def\u{{\bf u}}
\def\v{{\bf v}}
\def\l{\lambda}
\def\dsum{\displaystyle\sum}

\def\dsup{\displaystyle\sup}
\def\dint{\displaystyle\int}

\begin{document}

\title[quasineutral limit of  compressible Navier-Stokes-Poisson system]
{The quasineutral limit  of compressible Navier-Stokes-Poisson
system with heat
 conductivity and general initial data}

\author{Qiangchang Ju}
\address{Institute of Applied Physics and Computational Mathematics, PO
 Box 8009-28, Beijing 100088, P. R. China }
\email{qiangchang\_ju@yahoo.com}

 \author{Fucai Li}
\address{Department of Mathematics, Nanjing University,  Nanjing
 210093, P. R. China}
 \email{fli@nju.edu.cn}

\author{Hailiang li}
\address{Department of Mathematics and Institute of Mathematics and
Interdisciplinary Science, Capital Normal University, Beijing
100037, P. R. China} \email{hailiang.li.math@gmail.com}



\begin{abstract}
The quasineutral limit  of    compressible Navier-Stokes-Poisson
system with heat
 conductivity and general (ill-prepared) initial data
  is rigorously proved  in this paper. It is
proved that, as the   Debye length    tends  to zero, the solution
of the  compressible Navier-Stokes-Poisson system converges strongly
to the strong solution of the incompressible Navier-Stokes equations
plus a term of fast singular oscillating gradient vector fields.
Moreover,
  if the Debye length, the viscosity    coefficients  and the heat
conductivity coefficient independently go to zero, we obtain the
incompressible Euler equations.  In both cases  the convergence
rates are obtained.

\end{abstract}

\maketitle

\section{Introduction}
In the present paper we study the quasineutral limit of compressible
Navier-Stokes-Poisson system
\begin{align}
\partial_t \rho+\mbox{div}\,(\rho\u)&=0,\label{pns1}\\
\rho\{\partial_t\u+ (\u\cdot\nabla)\u\}+\nabla P(\rho,\theta)+\rho\nabla \Phi&=\mu\Delta \u
+(\mu+\nu)\nabla\mbox{div}\, \u , \label{ins2}\\
c_V\rho\{\partial_t \theta+(\u\cdot\nabla )\theta\}+P(\rho,\theta)
\mbox{div}\, \u&=\kappa \Delta \theta +\nu (\mbox{div}\,\u)^2+2\mu\mathbb{D}
(\u):\mathbb{D} (\u),\label{pns3}\\
 -\lambda^2\Delta\Phi  &=\rho -1, \label{pnsp}
\end{align}
  for $x\in \mathbb{T}^N\subset\mathbb{R}^N(N=2,3)$, the $N$-dimensional torus,  where $\rho , \u=(u_1,\dots,
  u_N),
$ $\theta$, and $\Phi $ denote the electron density, velocity,
temperature, and the electrostatic potential, respectively.
$\mathbb{D} (\u)=(d_{ij})_{i,j=1}^N, d_{ij}={1\over2}(\partial_i
u_{j}+\partial_j u_{ i})$.
 The constants $\nu$ and $\mu$ are the viscosity
coefficients with $\mu>0$ and $2\mu+N\nu>0$. $c_V>0$ is the specific
heat constant,
 $\kappa>0$ the heat conductivity coefficient, and $\lambda>0$ the scaled Debye length.
 The pressure function $P(\rho,\theta)$ takes the form
\begin{equation}\label{p}
P(\rho,\theta)=R\rho\theta, \quad R>0.
\end{equation}

 Without loss of generality, we assume $c_V=R\equiv 1$  for
notational simplicity. The Navier-Stokes-Poisson system
\eqref{pns1}-\eqref{pnsp} can be used to describe the dynamics of
plasma, where the compressible fluid of electron  interacts with its
own electric field against a charged ion background,
 see Degond~\cite{De00}.

    The purpose of the present paper is to investigate the
quasineutral limit of the  compressible  Navier-Stokes-Poisson
system \eqref{pns1}-\eqref{pnsp}. We shall prove rigorously that, as
the Debye length $\lambda\rightarrow 0$, the solution of the
compressible Navier-Stokes-Poisson system converges strongly to the
strong solution of the incompressible Navier-Stokes equations plus a
term of fast singular oscillating gradient vector fields   as long
as the strong solution of the latter   exists. Moreover, we also
consider the convergence of the compressible Navier-Stokes-Poisson
system \eqref{pns1}-\eqref{pnsp} to the incompressible Euler
equations
 by performing the combined quasineutral, vanishing
viscosity and vanishing  heat conductivity limit, i.e.
$\lambda\rightarrow 0$ and  $\mu, \nu,  \kappa\rightarrow 0$.

We first give   some  formal analysis. We use the subscript
$\lambda$ to indicate that the unknowns  are  dependent   on
$\lambda$ and set $\phi_\l= \l\Phi_\l$.   Thus, we can rewrite the
system \eqref{pns1}-\eqref{pnsp} as
 \begin{align}
&\partial_t\rho_\l+\mbox{div}(\rho_\lambda\u_\lambda)=0,\label{pns4}\\
& \rho_\lambda\{\partial_t\u_\lambda+ (\u_\l\cdot\nabla)\u_\l\}+
\nabla
(\rho_\lambda\theta_\lambda)+{1\over\l}\rho_\l\nabla\phi_\l=\mu\Delta\u_\l+(\nu+\mu)\nabla\mbox{div}
\u_\l,\label{pns5}\\
& \rho_\l\{\partial_t\theta_\l+ (\u_\l\cdot\nabla)\theta_\lambda\}
+\rho_\lambda\theta_\lambda\mbox{div}\u_\l\nonumber\\
&\qquad\qquad=
\kappa\Delta\theta_\l+\nu(\mbox{div}\u_\l)^2+2\mu\mathbb{D}
(\u_\l):\mathbb{D} (\u_\l),\label{pns6}\\
& -\lambda\Delta\phi_\l=\rho_\l-1.\label{pns7}
\end{align}

The system  \eqref{pns4}-\eqref{pns7} is equipped with the initial
data
\begin{equation} \label{ini}
\rho_\l(x,0)=\rho_{0\l}(x),\quad \u_\l(x,0)=\u_{0\l}(x),\quad
\theta_\l(x,0)=\theta_{0\l}(x).
\end{equation}

Letting $\l\rightarrow 0$ formally   in the Poisson equation
\eqref{pns7}, we have $\rho_\l=1$. Moreover, if we assume that
\begin{align}
\u_\l\rightarrow \v,\quad \theta_\l\rightarrow \theta  \nonumber
\end{align}
as $\l\rightarrow 0$, we may expect that the compressible
 Navier-Stokes-Poisson system
\eqref{pns4}-\eqref{pns7} converges to the incompressible
Navier-Stokes equations (see \cite{Lio96})
\begin{equation}\label{NS}
\left\{
\begin{aligned}
&\nabla\cdot \v=0,\\
&\partial_t\v+(\v\cdot\nabla)\v+\nabla\Pi=\mu\Delta\v,\\
&  \partial_t \theta+(\v\cdot\nabla)\theta =\kappa \Delta \theta  +
\displaystyle\frac{\mu}{2}\sum^{N}_{i,j=1}(\partial_iv_j+\partial_jv_i)^2,
\end{aligned}
\right.
\end{equation}
as the Debye length goes to zero, where $\nabla\Pi$ is expected to
be taken as the limit of the singular electric field and the
gradient of pressure  together. Furthermore, if  we let
$\mu\rightarrow 0$ and $\kappa \rightarrow 0$   in \eqref{NS}, it
yields the incompressible Euler equations
\begin{equation}\label{Euler}
\left\{
\begin{aligned}
&\nabla\cdot \v=0,\\
&\partial_t\v+(\v\cdot\nabla)\v+\nabla\Pi=0,\\
&  \partial_t \theta+(\v\cdot\nabla)\theta  =0.
\end{aligned}
\right.
\end{equation}

 Recently, there are  many progresses on the quasineutral limit of the compressible isentropic
  Navier-Stokes-Poisson system (i.e. the system \eqref{pns4}, \eqref{pns5} and \eqref{pns7}
  with the pressure $P_\lambda=a\rho^\gamma_\l, \gamma>1,a>0 $),
Wang~\cite{Wan04} studied the quasineutral limit for the smooth
solution with well-prepared initial data. Wang and
Jiang~\cite{WaJ06} studied the combined quasineutral and inviscid
limit of the compressible Navier-Stokes-Poisson system for weak
solution and obtained the convergence of Navier-Stokes-Poisson
system to the incompressible Euler equations with general initial
data. In\cite{WaJ06},  the vanishing of viscosity coefficients was
required in order to take the quasineutral   limit and no
convergence rate was derived therein. Ju, Li and Wang~\cite{JLyW}
improved the  arguments in \cite{WaJ06} and obtained the convergence
rate. Donatelli and Marcati~\cite{DM08} investigated the
quasineutral limit of the isentropic Navier-Stokes-Poisson system in
the whole space $\mathbb{R}^3$ and obtained the convergence of weak
solution of the Navier-Stokes-Poisson system to the weak solution of
the incompressible Navier-Stokes equations by means of dispersive
estimates of Strichartz's type under the assumption that the Mach
number is related to the Debye length. Notice that their arguments
can not be applied  to  the periodic case since the dispersive
phenomenon disappears in this situation. Ju, Li and
Wang~\cite{JLW08}  studied the quasineutral limit of the isentropic
Navier-Stokes-Poisson system both in the whole space and in the
torus  without restriction on  the viscosity coefficients.

 However, there is no  analysis on the quasineutral limit of the
  compressible non-isentropic Navier-Stokes-Poisson system  yet.
   In the present paper, we shall consider the  \emph{general
ill-prepared initial data} for the system \eqref{pns4}-\eqref{pns7},
 so  the fast oscillating singular term will be produced by the
non-divergence free part of initial momentum, and has to be
described carefully in order to pass  into the quasineutral limit.

In order to  describe the oscillations in time, we introduce the
following group $\mathcal{L}=e^{\tau L}, \tau \in \mathbb{R}$,
where $L$ is the operator defined on the space
$\mathcal{H}=(L^2(\mathbb{T}^N))^N\times\{\nabla\psi, \psi\in
H^1(\mathbb{T}^N)\}$ by
\begin{align}
L\left(
\begin{array}{c}
\bf{w}\\
0 \end{array} \right)=0,\;\;
\mbox{if}\;\;\mbox{div}\,\mathbf{w}=0,\nonumber\\
L\left(
\begin{array}{c}
\nabla q\\
\nabla\psi
\end{array}
\right)=\left(
\begin{array}{c}
-\nabla\psi\\
\nabla q
\end{array}
\right).
\end{align}
Then it is easy to check that $e^{\tau L}$ is an isometry on space
$H^s(\mathbb{T}^N)\times H^s(\mathbb{T}^N)$. Let us consider the
evolution of velocity and electric field. From \eqref{pns5} and
\eqref{pns7}, it is easy to obtain the following equation
\begin{equation}\label{ve}
\partial_t\nabla\phi_\lambda-{1\over\lambda}\mathcal{Q}\u_\lambda= -\mathcal{Q}(\u_\l\nabla\cdot
(\nabla\phi_\lambda)),
\end{equation}
  where the operator $\mathcal{Q}\v = \nabla \Delta^{-1}\nabla
\cdot \v$    is the Leray's projector on the space of gradient of
vector field $\v\in (L^2(\mathbb{T}^N))^N$, which  is defined as
follows
$$
 \mathcal{Q}\v = \nabla \Delta^{-1}\nabla
\cdot \v, \quad \mathcal{P}\v=(I-\mathcal{Q})\v,\quad \nabla\cdot
\mathcal{P}\v=0.
$$
We project the momentum equation \eqref{pns5} on the ``gradient
vector fields'' to obtain
\begin{align}\label{ve1}
 \partial_t\mathcal{Q}\u_\l+{1\over\lambda}\nabla\phi_\l=&-\mathcal{Q}((\u_\l\cdot\nabla)\u_\l)\nonumber\\
 & -\mathcal{Q}({1\over\rho_\l}\nabla
P_\l)+\mu\mathcal{Q}(\Delta\u_\l)+(\nu+\mu)\mathcal{Q}(\nabla\mbox{div}\u_\l)\nonumber\\
&+\mu\mathcal{Q}\Big(\big({1\over\rho_\l}-1)\Delta\u_\l\Big)+
(\nu+\mu)\mathcal{Q}\Big(\big({1\over\rho_\l}-1\big)\nabla\mbox{div}\u_\l\Big).
\end{align}
Define
$$U_\l=
\left(
\begin{array}{l}
\mathcal{\mathcal{Q}}\u_\l\\  \nabla\phi_\l
\end{array}
\right),\qquad V_\l=\mathcal{L}\Big(-\frac{t}{\l}\Big)U_\l.
$$
Then we can rewrite the system \eqref{ve}-\eqref{ve1} as
\begin{align}
\partial_tV_\l=\mathcal{L}\Big(-\frac{t}{\l}\Big)\left(
\begin{array}{l}
k_0\\   k_1
\end{array}\right),
\end{align}
with
\begin{align}
k_0=&-\mathcal{Q}((\u_\l\cdot\nabla)\u_\l)
-\mathcal{Q}({1\over\rho_\l}\nabla
P_\l)+\mu\mathcal{Q}(\Delta\u_\l)+(\nu+\mu)\mathcal{Q}(\nabla\mbox{div}\u_\l)\nonumber\\
&+\mu\mathcal{Q}\Big(\big({1\over\rho_\l}-1)\Delta\u_\l\Big)+
(\nu+\mu)\mathcal{Q}\Big(\big({1\over\rho_\l}-1\big)\nabla\mbox{div}\u_\l\Big),
 \label{k00}\\
k_1=&\mathcal{Q}\big(\u_\l\nabla\cdot
(\nabla\phi_\lambda)\big).\label{k11}
\end{align}

Now we can construct the oscillating terms as follows. Let $\v\in
C([0,T];$ $H^s(\mathbb{T}^N))$ be a divergence free function.
Consider the following linear system
\begin{align}\label{osc1}
\left\{ \begin{aligned}
 &\partial_t\nabla
q+\frac{1}{2}\mathcal{Q}\Big((\v\cdot\nabla)\nabla q+(\nabla q\cdot
\nabla)\v+\v\Delta
q\Big)-(\mu+\nu/2)\nabla\mbox{div}(\nabla q)=0,\\
&\partial_t\nabla p+\frac{1}{2}\mathcal{Q}\Big((\v\cdot\nabla)\nabla
p+(\nabla p\cdot \nabla)\v+\v\Delta
p\Big)-(\mu+\nu/2)\nabla\mbox{div}(\nabla p)=0
\end{aligned}
\right.
\end{align}
with initial data
$$(\nabla q(x,0), \nabla p(x,0))=(\mathcal{Q} {\u}_0(x),
\nabla\phi_0(x)).$$
It is direct to prove that there exists a unique
global smooth solution $(\nabla q,\nabla p)$ to the oscillating
system \eqref{osc1} satisfying
\begin{align}\label{qp}
\|(\nabla q,\nabla p)(t)\|_{H^s(\mathbb{T}^N)}\leq C(
T)\|(\mathcal{Q}{\u}_0,\nabla\phi_0)\|_{H^s(\mathbb{T}^N)},
\end{align}
where $C(T)>0$ is a constant depending only on $T$.

 Define
\begin{equation}\label{def1}
\left(
\begin{array}{l}
\u_{\rm osc}(x,t)\\
\nabla\phi_{\rm osc}(x,t)
\end{array}
\right)=\mathcal{L}(\frac{t}{\l})\left(
\begin{array}{l}
\nabla q(x,t)\\
\nabla p(x,t)
\end{array}
\right).
\end{equation}

 Before stating our results rigorously, we first recall the local
well-posedness result on the initial value problem for the
incompressible Navier-Stokes system \eqref{NS}   in multi-dimension.
One can refer to~\cite{Lio96} for the proof.

\begin{prop}\label{Proa}
Assume that $s\geq N/2+1$ and
\begin{equation}\label{initial3}
\left\{\begin{aligned}
&\v(x,0)=\v_0(x)\in H^{s+3}, \quad\ \mbox{\rm div}\, \v_0=0, \\
&\theta(x,0)=\theta_0(x)\in H^{s+3},\quad \ \inf_{x\in
\mathbb{T}^N}\theta_{0}(x)>0.
\end{aligned}\right.
\end{equation} Then there exists some time $T^*  ( 0<T^*\leq +\infty)$ such
that the initial problem   \eqref{NS} and \eqref{initial3} admits
a unique strong solution $(\v,\theta)$ satisfying,   for any
$T<T^*$,
\begin{align}
\v\in C^i([0,T], H^{s+3-i}), \quad i=0,1, \;\; \;\;\quad
\|\v(t)\|_{H^{s+3}}\leq
C_0\|\v_0\|_{H^{s+3}},\label{ves}\\
\theta\in C^i([0,T], H^{s+3-i}), \quad i=0,1, \;\; \;\;\quad
\|\theta(t)\|_{H^{s+3}}\leq C_0\|\v_0\|_{H^{s+3}}\label{ves1}
\end{align}
with $C_0>0$ a constant. Moreover, if $N=2$, the initial problem
\eqref{NS} and \eqref{initial3} admits a global unique strong
solution $(\v,\theta)\in C^i([0,\infty), H^{s+3-i}), i=0,1.$
\end{prop}

Our main results of this paper read as follows.

\begin{thm}\label{Thm}
Let $0< T<T^* $ defined in Proposition  \ref{Proa} and suppose that
 $(\v,\theta)\in C^i([0,T],$ $  H^{s+3-i})$, $ i=0,1,\; s>N/2+2
,$ be the unique strong solution of the initial problem  \eqref{NS}
and \eqref{initial3}. Assume that the initial data
$(\rho_{0\lambda}(x),\u_{0\lambda}(x), \theta_{0\lambda}(x))$
satisfies
\begin{align}
&\rho_{0\lambda}(x)=1-\lambda \Delta \phi_{0\l}(x), \quad \inf_{x\in
\mathbb{T}^N}\rho_{0\lambda}(x)>0, \quad \nabla \phi_{0\lambda}\in
H^{s+1}(\mathbb{T}^N),\label{iia}\\
& \qquad \u_{0\lambda}\in H^{s}(\mathbb{T}^N), \quad
\theta_{0\lambda}(x)\in H^{s}(\mathbb{T}^N), \quad \inf_{x\in
\mathbb{T}^N}\theta_{0\lambda}(x)>0,\label{iib}
 \end{align}
and
\begin{align}
&\qquad \|\mathcal{P}\u_{0\lambda}
-\v_0\|_{H^s}+\|\mathcal{Q}\u_{0\lambda}
-\mathcal{Q}\u_0\|_{H^s}\leq \tilde{C}\lambda\label{assin},\\
& \|\rho_{0\lambda}(x)-1  + \lambda \Delta \phi_0(x)\|_{H^s}\leq
\tilde{C}\lambda^2,  \quad \|\theta_{0\lambda}-\theta_0\|_{H^s}\leq
\tilde{C}\lambda\label{assin1}
\end{align}
for some constant $\tilde{C}>0$, where $\phi_0 $  and $  \u_0$ are
defined by \eqref{oscc}.
   Then there is a small constant $\delta_T>0$ such that, for any $
\lambda\in(0,\delta_T]$, the initial value problem for
Navier-Stokes-Poisson system \eqref{pns4}-\eqref{pns7}  admits a
unique classical solution
 $(\rho_\lambda, \u_\l, \theta_\l,\phi_\l)$ on $[0,T]$ satisfying
\begin{align} \sup_{0\leq t\leq T}\|(\rho_\lambda, \u_\l,
\theta_\l)(t)\|_{H^s}+\sup_{0\leq t\leq T}\|
\nabla\phi_\l(t)\|_{H^{s+1}}\leq C_1\label{uniform}
\end{align}
uniformly with respect to    $\lambda$.
 Moreover, it holds that
  \begin{align} & \sup_{0\leq t\leq
T}\big\{\|(\rho_\l-1)(t)\|_{H^s}+\|(\u_\l - \v-\u_{\rm
osc})(t)\|_{H^s} +\|(\theta_\l-\theta)(t)\|_{H^s}\big\}\nonumber\\
&\quad   + \sup_{0\leq t\leq T}\|(\nabla \phi_\l- \nabla\phi_{\rm
osc} )(t)\|_{H^{s+1}}  \leq C_2\lambda \label{limit1}
\end{align}
with $C_2>0$ independent of    $\lambda$.
\end{thm}

If we further perform  the  combined quasineutral, vanishing
viscosity and vanishing  heat conductivity limit, i.e.
$\lambda\rightarrow 0$  and $\mu,\nu, \kappa\rightarrow 0$, we
obtain
   the convergence of the
Navier-Stokes-Poisson system \eqref{pns1}-\eqref{pnsp} to the
incompressible Euler equations \eqref{Euler}. Namely,

 \begin{thm}\label{Thm2}
Let $0< T<T^{**} $ and suppose that
  $(\v,\theta)\in C^i([0,T],$ $  H^{s+3-i})$, $ i=0,1,\; s>N/2+2
,$ be the unique strong solution of the initial problem
\eqref{Euler} and \eqref{initial3}, where $T^{**}$ is the maximal
existing time of $(\v,\theta)$. Assume that the initial data
$(\rho_{0\lambda}(x),\u_{0\lambda}(x),\theta_{0\lambda}(x))$
satisfies the conditions \eqref{iia}-\eqref{assin1}.
  Then,    there is a small constant $\bar\delta_T>0$ such that, for any $
\lambda\in(0,\bar\delta_T]$, the initial value problem for
Navier-Stokes-Poisson system \eqref{pns4}-\eqref{pns7}  admits a
unique classical solution $(\rho_\lambda, \u_\l, \theta_\l,\phi_\l)$
on $[0,T]$ satisfying
\begin{align}\label{rate}
\sup_{0\leq t\leq T}\|(\rho_\lambda, \u_\l,
\theta_\l)(t)\|_{H^s}+\sup_{0\leq t\leq T}\|
\nabla\phi_\l(t)\|_{H^{s+1}}\leq C_3
\end{align}
uniformly with respect to    $\lambda$ as $\mu,\nu,\kappa\rightarrow
0$.
 Moreover, it holds that
  \begin{align} & \sup_{0\leq t\leq
T}\big\{\|(\rho_\l-1)(t)\|_{H^s}+\|(\u_\l - \v-\u_{\rm
osc})(t)\|_{H^s} +\|(\theta_\l-\theta)(t)\|_{H^s}\big\}\nonumber\\
&\quad   + \sup_{0\leq t\leq T}\|(\nabla \phi_\l- \nabla\phi_{\rm
osc} )(t)\|_{H^{s+1}}  \leq C_4\lambda \label{limit1111}
\end{align}
with $C_4>0$ independent of    $\lambda$.
 Here  $(\v,\theta)$ is  the unique strong solution of the initial problem
\eqref{Euler} and \eqref{initial3},  and  $(\u_{\rm osc},\phi_{\rm
osc})$ is the fast singular oscillating gradient velocity vector
field and electric field defined by \eqref{osc1} and \eqref{def1}
with $\mu=\nu\equiv0$.
\end{thm}

\begin{rem}
   The method developed in this paper can be applied to  the situation when the doping function is a perturbation of
a constant state
\begin{align}
\mathcal{C}(x)=1+\lambda g(x)\nonumber
\end{align}
with $g(x) \in C^2(\mathbb{T}^N)$, a given function, satisfying
$\int_{\mathbb{T}^N} g dx=0$.
\end{rem}

\begin{rem}
   We believe that the method developed in this paper can be also applied to investigate the quasineutral limit problem
   to  more complex model such  as the
  full   Navier-Stokes-Poisson system with more general pressure, which will be studied in a  forthcoming paper.
\end{rem}

The proofs of Theorems \ref{Thm} and \ref{Thm2} mainly consist of
three steps. First, we  apply the homogenization technique to
construct the approximate solution to the classical solution (if
exists) of the system \eqref{pns4}-\eqref{pns7}. Then by using the
theories of symmetric quasilinear hyperbolic system and the
estimates of second order elliptic equations, we show that the
remainder term exists in the same time interval as the approximate
term for fixed
 small $\lambda>0$. Moreover, we obtain the uniform estimates with
respect to $\lambda$  (the uniform estimates with respect to $\mu,
\nu$ and $\kappa$ can also be obtained by further analysis). These
facts are sufficient for us to complete the proofs of Theorems
\ref{Thm} and \ref{Thm2}.

It should be noted   that the quaineutral limit is a well-known
challenging and  modelling problem in fluid dynamics and kinetic
models for semiconductors and plasmas. In both cases there exist
only partial results. In particular, the quasineutral limit has been
performed in Vlasov-Poisson system by Brenier~\cite{Br},
Grenier~\cite{Gr1}, and Masmoudi~\cite{Ma}, in
Vlasov-Poisson-Fokker-Planck system by Hsiao, Li and
Wang~\cite{HLW06,HLW08}, in Schr\"odinger-Poisson system by
Puel~\cite{P1}, J\"ungel and Wang~\cite{JW}, and Ju et
al.~\cite{JLL}, in drift-diffusion-Poisson system    by Gasser et
al.~\cite{GLMS}, J\"ungel and Peng~\cite{JP}, Wang et al.
\cite{WXM}. For the hydrodynamic model, besides the results
mentioned above for the Navier-Stokes-Poisson system, there are also
many results on Euler-Poisson system,for example, for the isentropic
Euler-Poisson system~\cite{SS,CG1,Wan04,PW} and for non-isentropic
Euler-Poisson system~\cite{PWY,LLM}. Li and Lin~\cite{LL} considered
the quasineutral limit to  the isentropic quantum hydrodynamical
model with the help of modulated energy method
 for general initial data.

Before ending this section, we  recall the following Moser-type
calculus inequalities which will be used frequently in the sequel.

\begin{prop}[\!\!\cite{KM} Moser-type  inequalities]

(1) For $f,g\in H^s\cap L^\infty$ and $|\alpha|\leq s$, it holds
that
\begin{align}
\|D^\alpha(fg)\|_{L^2}\leq C_s(\|f\|_{L^\infty}\|D^s
g\|_{L^2}+\|g\|_{L^\infty}\|D^s\|_{L^2}).\label{MI1}
\end{align}

(2) For $f\in H^s, Df\in L^\infty, g\in H^{s-1}\cap L^\infty$ and
$|\alpha|\leq s$, it holds that
\begin{align}
\|D^\alpha(fg)-f D^\alpha(g)\|_{L^2}\leq
C_s(\|Df\|_{L^\infty}\|D^{s-1}g\|_{L^2}+\|g\|_{L^\infty}\|D^sf\|_{L^2}).\label{MI2}
\end{align}
\end{prop}

\emph{Notations}. In this paper, $C$ and $C_i(i=1,2,\dots)$ denote
the generic positive constants, which may change from line to line
and are independent of $\lambda$. $C(T)$ and $C_i(T)$ denote the
constant depending on the time $T$. $H^s$ denotes the standard
Sobolev space $W^{s,2}(\mathbb{T}^N)$. For the multi-index $\alpha =
(\alpha_1, \dots,\alpha_N)$,  we denote  $D^\alpha
=\partial^{\alpha_1}_{x_1}\cdots \partial^{\alpha_N}_{x_N}$ and
$|\alpha|=|\alpha_1|+\dots+|\alpha_N|$.

\smallskip
 The rest of this paper is arranged as follows. In Section
2, we construct  the approximate solutions to the problem
\eqref{pns4}-\eqref{ini}. In Section 3, we establish the local
existence of solution to the remainder system and obtain the uniform
estimates. The proofs of our main results are given in Section 4.

\section{Construction of approximate solutions}

In this section we shall construct the approximation to the system
\eqref{pns4}-\eqref{pns7}. Noticing the fast singular oscillating
vector fields  $(\u_{\rm osc}, \nabla\phi_{\rm osc})$ obtained by
\eqref{def1},   we find that the fast singular oscillating vector
fields $(\u_{\rm osc}, \nabla\phi_{\rm osc})$ satisfy
\begin{equation}\left\{
\begin{aligned}\label{oscc}
&\partial_t\u_{\rm
osc}+\frac{1}{2}\mathcal{Q}\Big((\v\cdot\nabla)\u_{\rm osc}+(\u_{\rm
osc}\cdot\nabla)\v+ \v\nabla\cdot\u_{\rm
osc}\Big)\\
&\qquad\qquad-(\mu+\nu/2)\nabla\mbox{div} \u_{\rm osc}
+\frac{1}{\lambda}\nabla\phi_{\rm osc}=0,\\
&\partial_t\nabla\phi_{\rm
osc}+\frac{1}{2}\mathcal{Q}\Big((\v\cdot\nabla)\nabla\phi_{\rm osc}+
(\nabla\phi_{\rm osc}\cdot\nabla)\v+ \v\Delta\phi_{\rm osc}\Big)\\
&\qquad\qquad-(\mu+\nu/2)\nabla\Delta\phi_{\rm osc}
-\frac{1}{\lambda}\u_{\rm osc}=0,\\
&(\u_{\rm osc}(x,0), \nabla\phi_{\rm osc}
(x,0))=(\mathcal{Q}{\u}_0(x), \nabla\phi_0(x)).
\end{aligned}
\right.
\end{equation}
Thus it is natural to define
\begin{align}
\rho_{\rm osc}=-\Delta \phi_{\rm osc}\nonumber.
\end{align}
We conclude that the fast oscillating part $(\rho_{\rm osc}, \u_{\rm
osc}, \phi_{\rm osc})$ satisfies the following initial value problem
\begin{equation}\label{osc11}
\left\{
\begin{aligned}
&\partial_t\rho_{\rm osc}+[\v+\u_{\rm osc}]\cdot\nabla\rho_{\rm osc}
+\frac{1}{\l}(1+\l\rho_{\rm osc})\nabla\cdot\u_{\rm osc}=k_2,\\
&\partial_t\u_{\rm osc}+([\v+\u_{\rm osc}]\cdot\nabla)\u_{\rm osc}
+(\u_{\rm osc}\cdot\nabla)\v+\frac{1}{\l}\nabla\phi_{\rm osc}=k_3,\\
&-\Delta\phi_{\rm osc}=\rho_{\rm osc},\\
&  \rho_{\rm osc}(x,0)=-\Delta \phi_0(x), \quad \u_{\rm
osc}(x,0)=\mathcal{Q}{\u}_0(x),
\end{aligned}
\right.
\end{equation}
where
\begin{align}
&k_2=\nabla\cdot(\rho_{\rm osc}[\v+\u_{\rm
osc}])+\frac{1}{2}\nabla\cdot\big((\v\cdot\nabla)\nabla\phi_{\rm osc}+
(\nabla\phi_{\rm osc}\cdot\nabla)\v+
\v\Delta\phi_{\rm osc}\big)\nonumber\\
&\qquad-(\mu+\nu/2)\Delta^2\phi_{\rm osc},\\
&k_3=\frac{1}{2}\mathcal{Q}\big((\v\cdot\nabla)\u_{\rm osc}+(\u_{\rm
osc}\cdot\nabla)\v-\v\nabla\cdot\u_{\rm osc}\big)
+(\u_{\rm osc}\cdot\nabla)\u_{\rm osc}\nonumber\\
&\qquad+\mathcal{P}\big((\v\cdot\nabla)\u_{\rm osc}+(\u_{\rm
osc}\cdot\nabla)\v\big)+(\mu+\nu/2)\nabla\mbox{div} \u_{\rm osc}.
\end{align}
Moreover, by  virtue of \eqref{qp} and \eqref{def1}, we obtain that
\begin{align}
\|k_2\|_{H^{s-2}(\mathbb{T}^N)}+\|k_3\|_{H^{s-2}(\mathbb{T}^N)}\leq
C\|( \nabla \phi_0, \mathcal{Q}{\u}_0,
\v_0)\|_{H^s(\mathbb{T}^N)},\label{k2}
\end{align}
 where the constant $C>0$ is
independent of $\lambda$. To approximate the classical solution
$W=(\rho_\l, \u_\l, \theta_\l, \phi_\l)^{\rm{T}}$ of the initial
value  problem \eqref{pns4}-\eqref{ini} for small $\l$, we still
need to introduce an additional correction term
\begin{align}
W_{\rm cor}=(\l\rho_{\rm cor}, \u_{\rm cor}, \theta_{\rm cor},
\phi_{\rm cor})^{\rm{T}}. \nonumber
\end{align}

By utilizing the fast singular oscillating part and the given
functions $k_2$ and $k_3$, we can construct $(\rho_{\rm cor},\u_{\rm
cor},\theta_{\rm cor}, \phi_{\rm cor})$ by solving the following
linear
 initial value
problem
\begin{equation}\label{corrr}
\left\{
\begin{aligned}
&\partial_\tau\u_{\rm cor}+\nabla\phi_{\rm cor}=k_4,\\
&\partial_\tau\nabla\phi_{\rm cor}-\u_{\rm cor}=\nabla(-\Delta)^{-1}k_2,\\
&
\rho_{\rm cor}=-\Delta \phi_{\rm cor},\\
&\partial_\tau \theta_{\rm cor}=k_5, \\
&(\u_{\rm cor}, \nabla\phi_{\rm cor}, \theta_{\rm
cor})(x,0)=(\mathbf{0}, \mathbf{0}, 0),
\end{aligned}
\right.
\end{equation}
where
\begin{align*}
  k_4=&-k_3-\nabla \theta+\mu \Delta\u_{\rm
osc}+(\mu+\nu)\nabla\mbox{div}\u_{\rm osc},\\
k_5=&-\u_{\rm osc}\cdot\nabla \theta-\theta\nabla\cdot\u_{\rm osc}+
\nu(\mbox{div} \u_{\rm osc})^2\\
& +\frac{\mu}{2}\sum^{N}_{i,j=1}(\partial_iv_j+\partial_jv_i+\partial_iu^j_{\rm osc}+\partial_ju^i_{\rm osc})^2.
\end{align*}
Here we recall  that $(\v,\theta)$ is the solution to the system
\eqref{NS}.

By virtue of \eqref{def1}, \eqref{ves}, \eqref{ves1} and
\eqref{k2}, it is easy to prove the following existence results of
solutions to the problems \eqref{osc11} and \eqref{corrr}.

\begin{prop}\label{Prob}
Let $T>0, T<T^*$ be given. Let $\v,\theta\in C^i([0,T],
H^{s+3-i}), i=0,1, s>1+N/2,$ be the solution to the initial value
problem \eqref{NS} and \eqref{initial3}. Then the problem
\eqref{osc11} admits a unique classical solution $(\rho_{\rm osc},
\u_{\rm osc}, \nabla\phi_{\rm osc})^{\rm T}$ for $t\in [0, T]$
satisfying
\begin{align}\label{oss} \|\rho_{\rm osc}(t)\|_{H^{s+2}}+ \|(\u_{\rm osc},
 \nabla\phi_{\rm
osc})(t)\|_{H^{s+3}(\mathbb{T}^N)}\leq C_T,
\end{align}
and the problem \eqref{corrr} admits a unique classical solution
$(\rho_{\rm cor}, \u_{\rm cor},\theta_{\rm cor}, \nabla\phi_{\rm
cor})^{\rm T} $ for $t\in [0, T]$ satisfying
\begin{align}\label{coo}
\|\rho_{\rm cor}(\tau)\|_{H^{s+1}}+ \|(\u_{\rm cor},\theta_{\rm
cor}, \nabla\phi_{\rm cor})(\tau)\|_{H^{s+2}(\mathbb{T}^N)}\leq
C_T,
\end{align}
 where $C_T>0$  depends only on $T$ and the initial data
$(\v_0,\theta_0, \mathcal{Q}{\u}_0, \nabla \phi_0)$, but is
independent of $\lambda$.
\end{prop}

  According to  Propositions \ref{Proa} and  \ref{Prob},   we
can make the following asymptotic expansions   of the solution
$(\rho_\l,\u_{\l},\theta_{\l},\phi_\l)$
\begin{equation}\label{expa}
\left\{
\begin{aligned}
&\rho_\l(x,t)=1+\l\rho_{\rm osc}(x,t)+\l^2(\Delta \Pi(x,t)+\rho_{\rm
cor}(x,t/\l))+\l^2\rho_{\rm rem}(x,t),\\
&\u_{\l}(x,t)=\v+\u_{\rm osc}(x,t)+\l\u_{\rm cor}(x,t/\l)+\l\u_{\rm rem}(x,t),\\
&\theta_{\l}(x,t)=\theta(x,t)+\l \theta_{\rm cor}(x,t/\l)+\l
\theta_{\rm rem}(x,t),\\
&\phi_\l(x,t)=\phi_{\rm osc}(x,t)+\l(\Pi(x,t)+\phi_{\rm
cor}(x,t/\l))+\l\phi_{\rm rem}(x,t).
\end{aligned}
\right.
\end{equation}
Substituting \eqref{expa} into the Navier-Stokes-Poisson system
\eqref{pns4}-\eqref{pns7}, using \eqref{NS}, \eqref{osc11} and
\eqref{corrr}, and by  tedious but direct  computations, we can
show that $(\rho_{\rm rem},\linebreak \u_{\rm rem}, \theta_{\rm
rem}, \phi_{\rm rem})$ solves the following initial value problem
\begin{equation}\label{rem1}
\left\{
\begin{aligned}
&\partial_t\rho_{\rm rem}+\u_\l\cdot\nabla\rho_{\rm rem}+\frac{1}{\l}{\rho_\l}\mbox{div}\u_{\rm rem}=h_0,\\
&\partial_t\u_{\rm rem}+(\u_\l\cdot\nabla)\u_{\rm rem}+\lambda
\frac{\theta_\lambda}{\rho_\lambda}\nabla\rho_{\rm rem}+\nabla
\theta_{\rm rem}\\
& \qquad\qquad-\mu\Delta\u_{\rm
rem}-(\mu+\nu)\nabla\mbox{div}\u_{\rm
rem}=-\frac{1}{\l}\nabla\phi_{\rm rem}
+\mathbf{f}_0,\\
&\partial_t \theta_{\rm rem}+\u_\l\cdot\nabla \theta_{\rm
rem}+\theta_\lambda\mbox{div}\u_{\rm rem}-\kappa\Delta  \theta_{\rm
rem}  =\l\nu(\mbox{div}\u_{\rm
rem})^2\\
&\qquad\qquad +\frac{\mu\l}{2}\sum^N_{i,j=1}(\partial_iu_{\rm
rem}^j+\partial_ju_{\rm rem}^i)^2+g_0,\\
&-\Delta \phi_{\rm rem}=\rho_{\rm rem}
\end{aligned}
\right.
\end{equation}
with initial data
\begin{equation}\label{ini2}
\left\{
\begin{aligned}
&\rho_{\rm
rem}(x,0)=\frac{1}{\lambda^2}\big[\rho_{0\lambda}(x)-1+\lambda\Delta
\phi_0(x)\big]-\Delta \Pi(x,0),\\
&\u_{\rm
rem}(x,0)=\frac{1}{\lambda}\big[\u_{0\lambda}(x)-\v_0(x)-\mathcal{Q}\u_0(x)\big],\\
&\theta_{\rm
rem}(x,0)=\frac{1}{\lambda}\big[\theta_{0\lambda}(x)-\theta_0(x)\big].
\end{aligned} \right.
\end{equation}
In \eqref{rem1}, we denote
\begin{align}
h_0=& -\u_{\rm rem}\cdot\nabla\rho_{\rm osc}-\rho_{\rm
rem}\nabla\cdot(\u_{\rm osc}+\l\u_{\rm cor})
-\nabla\cdot(\rho_{\rm osc}\u_{\rm cor})\nonumber\\
&-(\v+\u_{\rm osc}+\l\u_{\rm cor}+\l\u_{\rm rem})\cdot\nabla\rho_{\rm cor}
-\rho_{\rm cor}\nabla\cdot(\u_{\rm osc}+\l\u_{\rm cor})\nonumber\\
&-\Delta\Pi_t-(\nabla(\Delta\Pi))(\v+\u_{\rm osc}+\lambda\u_{\rm cor}+\l\u_{\rm rem})\nonumber\\
&-\Delta\Pi\mbox{div}(\u_{\rm osc}+\l\u_{\rm cor}),\label{defh0}\\
{\bf f}_0=&  {\bf f}_{01}+{\bf f}_{02},\label{deff0}\\
g_0=&g_{01}+g_{02}\label{defg0}
\end{align}
with
\begin{align}
{\bf f}_{01}=&-\big((\u_{\rm cor}+\u_{\rm
rem})\cdot\nabla\big)(\v+\u_{\rm osc})
-\big((\v+\u_{\rm osc}+\lambda\u_{\rm cor}+\lambda\u_{\rm rem})\cdot\nabla\big)\u_{\rm cor}\nonumber\\
&-   \frac{\theta_\lambda}{\rho_\lambda}   \nabla\big(\rho_{\rm
osc}+\lambda(\Delta\Pi+\rho_{\rm cor})\big)-\nabla
\theta_{\rm cor},\nonumber\\
{\bf f}_{02}=&\mu\Delta\u_{\rm cor}+(\mu+\nu)\nabla\mbox{div}\u_{\rm
cor}\nonumber\\
&  -\frac{\mu}{\rho_\l}\big(\rho_{\rm osc}
+\lambda(\Delta\Pi+\rho_{\rm cor})+\lambda\rho_{\rm
rem}\big)\Delta(\v+\u_{\rm osc}+\l\u_{\rm
cor}+\l\u_{\rm rem})\nonumber\\
&-\frac{\mu+\nu}{\rho_\l}\big(\rho_{\rm osc}
+\lambda(\Delta\Pi+\rho_{\rm cor})+\lambda\rho_{\rm
rem}\big)\nabla\mbox{div}(\v+\u_{\rm osc}+\l\u_{\rm
cor}+\l\u_{\rm rem}),\nonumber\\
g_{01}=&-(\u_{\rm cor}+\u_{\rm rem})\nabla \theta-(\v+\u_{\rm
osc}+\lambda\u_{\rm cor}+\lambda\u_{\rm rem})\nabla
\theta_{\rm cor}\nonumber\\
&-(\theta_{\rm cor}+\theta_{\rm rem})\mbox{div}\u_{\rm
osc}+\theta_\lambda\mbox{div}\u_{\rm cor},\nonumber\\
g_{02}=&\kappa\Delta\theta_{\rm cor}-\frac{\kappa}{\rho_\lambda}
\big(\rho_{\rm osc} +\lambda(\Delta\Pi+\rho_{\rm
cor})+\lambda\rho_{\rm rem}\big)\Delta(\theta+\lambda\theta_{\rm
cor}+\lambda\theta_{\rm
rem})\nonumber\\
& +2\nu\mbox{div}\u_{\rm osc}(\mbox{div}\u_{\rm
cor}+\mbox{div}\u_{\rm rem}) +\l\nu(\mbox{div}\u_{\rm cor})^2
+2\l\nu\mbox{div}\u_{\rm cor}\mbox{div}\u_{\rm rem}\nonumber\\
& +\mu\sum^N_{i,j=1}(\partial_iv_j+\partial_jv_i+\partial_iu_{\rm
osc}^j+\partial_ju_{\rm osc}^i) (\partial_iu_{\rm cor}^j+
\partial_ju_{\rm cor}^i+\partial_iu_{\rm rem}^j+\partial_ju_{\rm rem}^i)
\nonumber\\
& +\frac{\mu\l}{2}\sum^N_{i,j=1}(\partial_iu_{\rm cor}^j+
\partial_ju_{\rm cor}^i)^2+\mu\l\sum^N_{i,j=1}(\partial_iu_{\rm cor}^j+
\partial_ju_{\rm cor}^i)(\partial_iu_{\rm rem}^j+\partial_ju_{\rm rem}^i)\nonumber\\
& -\frac{1}{\rho_\lambda}\big(\rho_{\rm osc}+\lambda \Delta \Pi
+\lambda \rho_{\rm cor}+\lambda\rho_{\rm
rem}\big)\bigg[\nu(\mbox{div}\u_\lambda)^2
+\frac{\mu}{2}\sum^N_{i,j=1}(\partial_iu_\lambda^j+\partial_ju_\lambda^i)^2\bigg]
.\nonumber
\end{align}

If we denote
\begin{align*}
U_{\rm rem}:=(\rho_{\rm rem}, \u_{\rm rem}, \theta_{\rm rem})^{\rm
T},
\end{align*}
the problem \eqref{rem1}-\eqref{ini2} can be rewritten as follows
\begin{equation}\label{rem2}
\left\{
\begin{aligned}
&\partial_tU_{\rm rem}+\dsum_{j=1}^N A_j(x,t,U_{\rm
rem})\partial_{x_j}U_{\rm rem}-\mu\Delta\tilde{\u}_{\rm rem}
-(\mu+\nu)\nabla\mbox{div}\tilde{\u}_{\rm rem}\\
&\qquad -\kappa\Delta\tilde{\theta}_{\rm rem}=\lambda\nu J+\frac{\l\mu}{2} G+\dfrac{1}{\l}B+F(x,t, U_{\rm rem}),\\
&-\Delta \phi_{\rm rem}=\rho_{\rm rem}, \\
&U_{\rm rem}(x,0)=(\rho_{\rm rem}(x,0), \u_{\rm
rem}(x,0),\theta_{\rm rem}(x,0) )^{\rm T} := U_{\rm{rem}0}(x).
\end{aligned} \right.
\end{equation}

Here the matrices $A_j(j=1,\dots,N)$ is defined as
$$
A_j(x,t, U_{\rm rem})\equiv u_\l^j I_{(N+2)\times(N+2)}+ \left(
\begin{array}{ccc}
0&\dfrac{1}{\l}\rho_\l e_j,&0\\
 \frac{\l\theta_\lambda}{\rho_\lambda}e_j^{\rm T}& O& e_j^{\rm T}\\
0&\theta_\lambda e_j& 0
\end{array}
\right)$$
 and
 \begin{align*}
&\tilde{\u}_{\rm rem}  =(0, \u_{\rm rem}, 0)^{\rm{T}},
 &&\tilde{\theta}_{\rm rem}=(0,\dots,0,\theta_{\rm rem})^{\rm{T}},\\
&J=  (0,\dots, 0, (\mbox{div}\u_{\rm rem})^2)^{\rm{T}},    & &F =(h_0,\mathbf{f}_0,g_0)^{\rm{T}}, \\
&G= \bigg(0,\dots,0,\sum^N_{i,j=1}(\partial_iu_{\rm rem}^j+\partial_ju_{\rm
rem}^i)^2\bigg)^{\rm{T}},    &&B =(0,-\nabla\phi_{\rm rem}, 0)^{\rm{T}}.
\end{align*}

\section{Local existence of solution to the remainder system \eqref{rem2}}

In this section we study the local existence of smooth solution to
the remainder system \eqref{rem2}, our result reads

\begin{thm}\label{Thm3}
Let $T>0,T<T^*$ be given and $\v,\theta\in C^i([0,T], H^{s+3-i}),
i=0,1, s>2+N/2,$ be the solution to the problem \eqref{NS} and
\eqref{initial3}. Then there exists a  constant $\delta_T>0$ such
that for any $ \lambda\in(0,\delta_T]$, the initial value problem
\eqref{rem2}  admits a unique classical solution $(U_{\rm rem},
\phi_{\rm rem})$ in $[0,T]$ satisfying
\begin{align}
\dsup_{0\leq t\leq T}\big(\|(\lambda\rho_{\rm rem}, \u_{\rm rem},
\theta_{\rm rem})(t)\|_{H^s}+\|\nabla\phi_{\rm
rem}(t)\|_{H^{s+1}}\big)\leq C(T) ,\label{apriori2}
\end{align}
where $C(T)$ is a positive constant independent of $\lambda$.
\end{thm}

The proof of Theorem  \ref{Thm3}  proceeds via    a priori  energy
estimates and  the   classical iteration scheme. The crucial step
is to show the following energy estimates which can be obtained by
performing the refined  energy estimates for the quasilinear
symmetric hyperbolic-parabolic system and the Poisson equation.

\begin{lem}\label{apriori}
Let $T>0$ be given and $s\geq N/2+2$. There exist positive constants
$\delta_T, M, \tilde{M}$ such that the classical solutions $(U_{\rm
rem}, \phi_{\rm rem})$ to the initial value problem \eqref{rem2}
satisfies
\begin{align} &\dsup_{0\leq t\leq
T}\big(\|(\lambda\rho_{\rm rem}, \u_{\rm rem},
\theta_{\rm rem})(t)\|^2_{H^s}+\|\nabla\phi_{\rm rem}(t)\|^2_{H^{s+1}}\big)\nonumber\\
&\quad+ \dint_0^T\|\u_{\rm rem}(s)\|_{H^{s+1}}^2\;dt +
\dint_0^T\|\theta_{\rm rem}(s)\|_{H^{s+1}}^2\;dt   \leq M^2,
\label{apriori1}\end{align}
and
\begin{align}
&\dsup_{0\leq t\leq T}\Big(\|\lambda\partial_t\rho_{\rm
rem}(t)\|_{H^{s-1}}+\|\lambda\partial_t\u_{\rm rem}(t)\|_{H^{s-2}}
+\|\partial_t \theta_{\rm
rem}(t)\|_{H^{s-2}}\nonumber\\
&\qquad\qquad +\|\lambda\partial_t\nabla\phi_{\rm
rem}(t)\|_{H^s}\Big)\leq \tilde{M}\label{apriorit}
\end{align}
uniformly with respect to $\lambda\in (0, \delta_T]$.
\end{lem}

\begin{proof}[ Proof of Lemma \ref{apriori}] We assume a priori that the
classical solution to initial value problem \eqref{rem2} satisfies
\eqref{apriori1} and \eqref{apriorit}. Then our task is to determine
these unknown constants by energy estimates.

Noticing the matrices  $ A_j(x, t, U_{\rm rem}), j=1,\dots,N$ can be
symmetrized by
$$A_0(x,t,U_{\rm rem})=
\left(
\begin{array}{ccc}
\l^2\frac{\theta_\lambda}{\rho_\lambda}& O & 0\\
O  &\rho_\l I_{N\times N}& O\\
0& O & \frac{\rho_\l}{\theta_\lambda}
\end{array}
\right),
$$
we rewrite the system \eqref{rem2} in the following form
\begin{equation}\label{rem3}
\left\{
\begin{aligned}
&A_0(U_{\rm rem})\partial_tU_{\rm
rem}+\dsum_{j=1}^N\mathcal{A}_j(x,t,U_{\rm rem})\partial_{x_j}U_{\rm
rem}-
\mu\rho_\l\Delta\tilde{\u}_{\rm rem}\\
&\qquad  -(\mu+\nu)\rho_\l\nabla\mbox{div}\tilde{\u}_{\rm rem}
-\frac{\kappa\rho_\lambda}{\theta_\lambda}\Delta\tilde{\theta}_{\rm rem}\\
&\qquad\qquad=
 \lambda\nu \tilde{J}+\frac{\l\mu}{2} \tilde{G}+\dfrac{1}{\l}\tilde{B}+\tilde{F}(x,t, U_{\rm rem}),\\
&-\Delta \phi_{\rm rem}=\rho_{\rm rem}, \\
& U_{\rm rem}(x,0)=U_{{\rm rem}0}(x),
\end{aligned} \right.
\end{equation}
where $\mathcal{A}_j=A_0A_j, j=1,\dots,N$ are symmetric matrices
given by
$$\mathcal{A}_j(U_{\rm rem})= u_\l^j  A_0(U_{\rm rem})+
\left(
\begin{array}{ccc}
0&\l\theta_\lambda e_j& 0\\
\l\theta_\lambda e_j^{\rm T}& O& \rho_\l e_j^{\rm T}\\
0&\rho_\l e_j&0
\end{array}
\right)$$ and
\begin{align*}
&\tilde{J}:=A_0J=\Big(0,\dots,0,\frac{\rho_\lambda}{\theta_\lambda}(\mbox{div}\u_{\rm rem})^2\Big)^{\rm{T}},\\
&\tilde{G}:=A_0G=\bigg(0,\dots,0,\frac{\rho_\lambda}{\theta_\lambda}
\sum^N_{i,j=1}(\partial_iu_{\rm rem}^j+\partial_ju_{\rm rem}^i)^2\bigg)^{\rm{T}},\\
&\tilde{B}:=A_0B=(0, -\rho_\l\nabla\phi_{\rm rem}, 0)^{\rm{T}},\\
&\tilde{F}:=A_0F=\Big(\frac{\lambda^2\theta_\l h_0}{\rho_\l},
\rho_\l \mathbf{f}_0,\frac{\rho_\l g_0}{\theta_\l}\Big)^{\rm{T}}.
\end{align*}

Next we perform energy estimates for the classical solution to the
system \eqref{rem2} with initial data \eqref{ini2}. Define the
canonical energy by
\begin{align*}
\|U_{\rm rem}\|_{E}^2:=\dint \langle A_0(U_{\rm rem})U_{\rm rem},
U_{\rm rem}\rangle dx.
\end{align*}
Multiplying \eqref{rem3}$_1$ by $U_{\rm rem}$ and integrating the
result by parts, we get the basic energy equality of Friedrich's
\begin{align}
&\dfrac{d}{dt}\|U_{\rm rem}\|_{E}^2+2\mu\dint|\nabla\u_{\rm
rem}|^2dx+ 2(\mu+\nu)\dint|\mbox{div}\u_{\rm rem}|^2
dx\nonumber\\
&\quad+2\kappa\dint\frac{\rho_\lambda}{\theta_\lambda} |\nabla\theta_{\rm rem}|^2 dx\nonumber\\
=& \dint\langle \Gamma U_{\rm rem}, U_{\rm rem} \rangle dx
+2\lambda\nu \dint \frac{1}{\theta_\lambda}(\mbox{div}
\u_{\rm rem})^2\theta_{\rm rem}dx\nonumber\\
&+\lambda\mu\sum^N_{i,j=1} \dint
\frac{1}{\theta_\lambda}(\partial_i u_{\rm rem}^j+\partial_j
u_{\rm rem}^i)^2\theta_{\rm rem} dx
-\dfrac{2}{\lambda}\dint\rho_{\lambda}\nabla\phi_{\rm rem}\u_{\rm
rem}
dx\nonumber\\
&+2\dint \langle A_0F, U_{\rm rem}\rangle dx+ R_1,  \label{basic}
\end{align}
where
\begin{align}
R_1=& 2(\mu+\nu)\dint (\rho_\lambda-1)\nabla\mbox{div}\u_{\rm
rem}\u_{\rm rem}dx+2\mu\dint
(\rho_\lambda-1)\Delta\u_{\rm rem} \u_{\rm rem} dx\nonumber\\
&-2\kappa\dint
\nabla\big(\frac{\rho_\lambda}{\theta_\lambda}\big)\nabla\theta_{\rm
rem}\theta_{\rm rem} dx
\end{align}
and
\begin{align*}
\Gamma=(\partial_t, \nabla)\cdot(A_0,\mathcal{A}_1,\dots,
 \mathcal{A}_3).
\end{align*}

Since $\mu>0, 2\mu+N\nu>0$, there exists  a positive constant $\xi_1$
such that
\begin{align}
\mu\dint|\nabla\u_{\rm rem}|^2dx+ (\mu+\nu)\dint|\mbox{div}\u_{\rm
rem}|^2 dx\geq \xi_1 \dint |\nabla\u_{\rm rem}|^2dx\label{xi}
\end{align}
in view of $\int(\mbox{div} \u_{\rm rem})^2dx\leq \int
|\nabla\u_{\rm rem}|^2 dx.$ Notice the fact that there is a
$\delta_T>0$ such that for $\lambda\in(0, \l_T]$ it holds that
\begin{align}
&0<\rho_-\leq 1+\l\rho_{\rm osc}+\l^2\Delta\Pi+\l^2\rho_{\rm
cor}+\l^2\rho_{\rm rem}\leq
\rho_+,\label{rho1}\\
&\qquad\qquad 0<\theta_-\leq \theta+\lambda \theta_{\rm cor}+\l
\theta_{\rm rem}\leq \theta_+, \label{e1}
\end{align}
where $\rho_{\pm}$ and $\theta_{\pm}$ are positive constants.
Thus, the matrices $A_0$ and $\mathcal{A}_j, j=1,\dots,N,$
together with their derivatives are continuous and bounded
uniformly. Moreover, $A_0$ is uniformly positive definite, i.e.
there exists a constant $c_0>0$ such that
\begin{align}
\label{positive1} \langle A_0(U_{\rm rem})U_{\rm rem}, U_{\rm
rem}\rangle \geq c_0(\l^2\rho_{\rm rem}^2+\u_{\rm rem}^2+\theta_{\rm
rem}^2)
\end{align}
for all $U_{\rm rem}$.

Now we estimate the terms on the right-hand side of \eqref{basic}.
Since $\Gamma$ is bounded there exists a generic constant $M_0$,
  independent of $(\rho_{\rm rem}, \u_{\rm rem}, \theta_{\rm
rem}, \phi_{\rm rem})$ and $\l>0$, such that
\begin{align}
\dint\langle \Gamma U_{\rm rem},  U_{\rm rem}\rangle dx \leq
M_0(1+\lambda(M+\tilde{M}))\| U_{\rm rem}\|_{E}^2.\label{Gamma0}
\end{align}
By Sobolev's embedding inequality and the inequality \eqref{e1} we
obtain that
\begin{align}
&2\lambda\nu\dint \frac{1}{\theta_\lambda}(\mbox{div} \u_{\rm
rem})^2\theta_{\rm rem}dx +\lambda\mu\sum^N_{i,j=1}
\dint\frac{1}{\theta_\lambda}(\partial_i u_{\rm rem}^j+\partial_j
u_{\rm rem}^i)^2\theta_{\rm rem} dx\nonumber\\
\leq&\lambda M_0M(2\mu+\nu)\int (|\nabla\u_{\rm rem}|^2+
|\theta_{\rm rem}|^2)dx. \label{non1}
\end{align}

By integrating by parts, Cauchy's inequality and the equation for
$\rho_{\rm rem}$ in \eqref{rem2}, the forth term on the right-hand
side of \eqref{basic} is estimated as follows
\begin{align}
&-\dfrac{2}{\lambda}\dint\rho_{\lambda}\nabla\phi_{\rm rem}\u_{\rm
rem}
dx\nonumber\\
=&\frac{2}{\lambda}\dint \rho_\lambda\mbox{div}\u_{\rm rem}\phi_{\rm
rem}dx+\frac{2}{\lambda}\dint
\nabla\rho_\lambda\u_{\rm rem}\phi_{\rm rem} dx\nonumber\\
=&-2\dint\partial_t\rho_{\rm rem}\phi_{\rm rem} dx-2\dint(\v+\u_{\rm
osc}+\lambda\u_{\rm cor}+\lambda\u_{\rm rem})\nabla\rho_{\rm
rem}\phi_{\rm rem}
dx\nonumber\\
&+2\dint\dint h_0\phi_{\rm rem}dx+2\dint \nabla  (\rho_{\rm
osc}+\lambda(\Delta
\Pi+\rho_{cor})+\lambda\rho_{\rm rem})\u_{\rm rem}\phi_{\rm rem} dx\nonumber\\
\leq & -\partial_t\|\nabla\phi_{\rm rem}\|_{L^2}^2+M_0(1+\lambda
M)(\|\nabla\phi_{\rm rem}\|_{L^2}^2+\|U_{\rm
rem}\|_{E}^2)+\epsilon_1 \|\nabla\u_{\rm rem}\|_{L^2}^2
\label{forth}  \end{align} for some sufficiently small constant
$\epsilon_1>0$.

Now we deal with the term $R_1$. By integrating by parts and  using
Sobolev's inequality, we get
\begin{align}
&2(\mu+\nu)\dint (\rho_\lambda-1)\nabla\mbox{div}\u_{\rm rem}\u_{\rm
rem}dx+2\mu\dint
(\rho_\lambda-1)\Delta\u_{\rm rem} \u_{\rm rem} dx\nonumber\\
\leq&\lambda M_0(M+1)(2\mu+\nu)\int (|\nabla\u_{\rm rem}|^2+|\u_{\rm
rem}|^2)dx.\label{mu1}
\end{align}
In view of \eqref{rho1}, \eqref{e1} and Cauchy's inequality, we
obtain that
\begin{align}
&-2\kappa\dint
\nabla\frac{\rho_\lambda}{\theta_\lambda}\nabla\theta_{\rm
rem}\theta_{\rm rem}
dx\nonumber\\
=&-2\kappa\dint
\frac{\nabla\rho_\lambda}{\theta_\lambda}\nabla\theta_{\rm
rem}\theta_{\rm rem} dx+2\kappa\dint
\frac{\rho_\lambda}{(\theta_\lambda)^2}\nabla\theta_\lambda\nabla\theta_{\rm
rem}\theta_{\rm rem}
dx\nonumber\\
\leq & \lambda M_0(M+1)\kappa\dint (|\nabla\theta_{\rm
rem}|^2+|\theta_{\rm rem}|^2) dx \nonumber\\
& +M_0\kappa\dint |\theta_{\rm
rem}|^2 dx+\epsilon_2\kappa \dint |\nabla\theta_{\rm rem}|^2
dx\label{thetarem}
\end{align}
for some sufficiently small constant $\epsilon_2>0$.

The estimate of the fifth term on the right-hand side of
\eqref{basic} is  tedious but straightforward. In view of the
definitions of $h_0, {\bf f}_0,$ and $g_0$ in
\eqref{defh0}-\eqref{defg0}, and the Propositions \ref{Proa} and
\ref{Prob}, we get
\begin{align}
&2\lambda^2\dint \frac{\theta_\lambda}{\rho_\lambda}h_0\rho_{\rm
rem} dx +2\dint \rho_\lambda{\bf f}_{01}\u_{\rm rem} dx+2\dint
\frac{\rho_\lambda}{\theta_\lambda}
 g_{01} \theta_{\rm rem}dx\nonumber\\
\leq & M_0\|U_{\rm rem}\|_{E}^2+M_0 \label{01}
\end{align}
and
\begin{align}
&2\dint \rho_\lambda{\bf f}_{02}\cdot\u_{\rm rem}dx+2\dint
\frac{\rho_\lambda}{\theta_\lambda}
  g_{02} \theta_{\rm rem}     dx\nonumber\\
\leq &
\lambda (2\mu+\nu+\kappa)M_0(1+M)\dint (|\u_{\rm rem}|^2\nonumber
+|\nabla\u_{\rm rem}|^2+|\nabla\theta_{\rm rem}|^2)dx \\
&
+(2\mu+\nu+k+1)M_0.\label{F1}
\end{align}

We choose $\delta_T$ sufficiently  small such that, for
$\lambda\in (0, \delta_T]$,
\begin{align}
\lambda M_0(M+1)(2\mu+\nu+\kappa)\leq \min\Big\{\frac{\xi_1}{2},
\frac{\kappa\rho_-}{2\theta_+}\Big\}:=\eta_1. \label{lamsmall0}
\end{align}
Choosing $\epsilon_1 $ and $\epsilon_2 $ sufficiently small  and
combining \eqref{xi}-\eqref{lamsmall0} with \eqref{basic}, we
obtain that
\begin{align}
&\dfrac{d}{dt}\Big(\|U_{\rm rem}\|_{E}^2+\|\nabla\phi_{\rm
rem}\|_{L^2}^2 \Big)+\frac{\xi_1}{2}\dint  |\nabla\u_{\rm rem}|^2
dx+\frac{k\rho_-}{2\theta_+}\dint|\nabla\theta_{\rm rem}|^2
dx\nonumber\\
\leq & M_0(1+\lambda(M+\tilde{M}))(\| U_{\rm
rem}\|_{E}^2+\|\nabla\phi_{\rm rem}\|_{L^2}^2)+3\eta\dint(|\u_{\rm
rem}|^2 +|\theta_{\rm rem}|^2) dx\nonumber\\&+\kappa M_0\dint
|\theta_{\rm rem}|^2 dx +(2\mu+\nu+\kappa+1)M_0.\label{estimatezero}
\end{align}

Next we shall obtain the energy estimates of higher order
derivatives for the classical solutions to the initial value problem
\eqref{rem2}. For the multi-index $\alpha$ with $1\leq |\alpha|\leq
s$, we take the operator $D^\alpha$ to \eqref{rem2} and multiply the
resulting equations by $A_0$ to obtain
\begin{equation}\label{rem4}
\left\{
\begin{aligned}
&A_0(U_{\rm rem})\partial_tD^\alpha U_{\rm
rem}+\dsum_{j=1}^N\mathcal{A}_j(x,t,U_{\rm
rem})\partial_{x_j}D^\alpha
U_{\rm rem}-\rho_\l\mu\Delta D^\alpha\tilde{\u}_{\rm rem}\\
&\qquad-(\mu+\nu)\rho_\l\nabla\mbox{div}D^\alpha\tilde{\u}_{\rm rem}
-\frac{\kappa\rho_\lambda}{\theta_\lambda}\Delta D^\alpha\tilde{\theta}_{\rm rem}\\
&\quad = \lambda\nu A_0(U_{\rm rem}) D^\alpha J+\frac{\l\mu}{2}
A_0(U_{\rm rem}) D^\alpha G +\dfrac{1}{\l}A_0(U_{\rm rem}) D^\alpha B\\
&\qquad +A_0(U_{\rm rem})D^\alpha F
+ H^\alpha,\\
&-\Delta D^\alpha\phi_{\rm rem}=D^\alpha\rho_{\rm rem}
\end{aligned} \right.
\end{equation}
with initial data
\begin{equation}\label{ini4}
  D^\alpha U_{\rm rem}(x,0)=D^\alpha U_{\rm{rem}0}(x),
  \end{equation}
   where $H^\alpha $ consists of the commutating terms as
\begin{align*}
&H^\alpha =-\dsum_{j=1}^N A_0(U_{\rm rem})\Big(D^\alpha(A_j(U_{\rm
rem})\partial_{x_j}U_{\rm rem}) -A_j(U_{\rm
rem})\partial_{x_j}D^\alpha U_{\rm rem}\Big).
\end{align*}

Taking the inner product between \eqref{rem4}$_1$ and $D^\alpha U_{\rm
rem}$, we have the following differential equality
\begin{align}
&\dfrac{d}{dt}\|D^\alpha U_{\rm rem}(t)\|_{E}^2+2\mu\dint | \nabla
D^{\alpha} \u_{\rm rem}|^2 dx+2(\mu+\nu)\dint|\mbox{div} D^{\alpha}
\u_{\rm rem}|^2dx\nonumber\\
&\qquad+2\kappa\dint\frac{\rho_\lambda}{\theta_\lambda}\Big| D^{\alpha+1}\theta_{\rm rem}\Big|^2dx\nonumber\\
=&\dint\langle \Gamma D^\alpha U_{\rm rem}, D^\alpha U_{\rm
rem}\rangle dx +2\l\nu\dint\langle A_0(U_{\rm rem})D^\alpha J,
D^\alpha U_{\rm rem}\rangle
dx\nonumber\\
&+\l\mu\dint\langle A_0(U_{\rm rem})D^\alpha G,D^\alpha U_{\rm
rem}\rangle dx +\dfrac{2}{\l}\dint \langle A_0(U_{\rm rem}) D^\alpha
B,D^\alpha U_{\rm rem} \rangle dx\nonumber\\&+2\dint\langle
A_0(U_{\rm rem})D^\alpha F(t), D^\alpha U_{\rm rem}\rangle dx
+2\dint\langle H^\alpha(t), D^\alpha U_{\rm rem} \rangle dx+  R_2,
\label{higher1}
\end{align}
where
\begin{align}
R_2 = &2\mu\dint(\rho_\l-1)\Delta D^\alpha\u_{\rm
rem}D^\alpha\u_{\rm rem}
dx-2\kappa\dint\nabla\big(\frac{\rho_\lambda}{\theta_\lambda}\big)\nabla
D^\alpha\theta_{\rm rem}\theta_{\rm rem}dx \nonumber\\
&+2(\mu+\nu)\dint(\rho_\l-1)\nabla\mbox{div}D^\alpha \u_{\rm rem}
D^\alpha \u_{\rm rem}  dx.\nonumber
\end{align}

It is easy to see that we also have the following estimate
\begin{align}
\mu\dint|\nabla D^\alpha\u_{\rm rem}|^2dx+
(\mu+\nu)\dint|\mbox{div}D^\alpha\u_{\rm rem}|^2 dx\geq \xi_2 \dint
|\nabla D^\alpha\u_{\rm rem}|^2dx\label{xial}
\end{align}
for some constant $\xi_2>0$.

Now we deal with the right-hand side of \eqref{higher1}. In the
following the generic constant $M_0$ may depend on $T$ and $s$. By
 integrating by part, Sobolev's inequality and Cauchy's inequality
it holds, similar to \eqref{Gamma0} and
\eqref{mu1}-\eqref{thetarem}, that
\begin{align}
\dint\langle \Gamma D^\alpha U_{\rm rem}, D^\alpha U_{\rm
rem}\rangle dx \leq M_0(1+\lambda(M+\tilde{M}))\|D^\alpha U_{\rm
rem}\|_{E}^2\label{Gamma}
\end{align}
and
\begin{align}
R_2\leq &\lambda M_0(M+1)(2\mu+\nu+\kappa)\int \Big(  |\nabla
D^{\alpha}\u_{\rm rem}|^2
+|D^{\alpha+1}\theta_{\rm rem}|^2+|D^\alpha\u_{\rm rem}|^2\nonumber\\
&+|D^\alpha\theta_{\rm rem}|^2\Big)dx+M_0\kappa\dint
|D^\alpha\theta_{\rm rem}|^2 dx+\delta\kappa \dint
|D^{\alpha+1}\theta_{\rm rem}|^2 \label{I1}
\end{align}
for some sufficiently small constant $\delta >0$.

By the definition of $A_0,G$ and $J$, it follows from the
Sobolev's inequality that
\begin{align}
&\l\nu\dint\langle A_0(U_{\rm rem})D^\alpha J, D^\alpha U_{\rm
rem}\rangle dx+2\l\mu\dint\langle A_0(U_{\rm rem})D^\alpha
G,D^\alpha
U_{\rm rem}\rangle dx\nonumber\\
\leq & \lambda
M_0(2\mu+\nu)\bigg(\|({\rm div}\u_{\rm rem})^2\|_{H^\alpha}
+\Big\|\sum^N_{i,j=1}(\partial_iu_{\rm{rem}}^j
+\partial_ju_{\rm rem}^i)^2\Big\|_{H^\alpha}\bigg)\|\theta_{\rm rem}\|_{H^\alpha}\nonumber\\
\leq & \lambda M_0(2\mu+\nu)\|\u_{\rm
rem}\|_{H^\alpha}^{\frac{1}{4}}
\|\u_{\rm rem}\|_{H^{\alpha+1}}^{\frac{7}{4}}\|\theta_{\rm rem}\|_{H^\alpha}\nonumber\\
\leq &\lambda M_0(2\mu+\nu)\|\u_{\rm rem}\|_{H^\alpha}^{\frac{1}{2}}
\|\u_{\rm rem}\|_{H^{\alpha+1}}^{\frac{3}{2}}\|\theta_{\rm
rem}\|_{H^\alpha}^2+\l
M_0(2\mu+\nu)\|\u_{\rm rem}\|_{H^{\alpha+1}}^2\nonumber\\
\leq & \lambda M M_0(2\mu+\nu)\|\u_{\rm
rem}\|_{H^{\alpha+1}}^{\frac{3}{2}}\|\theta_{\rm
rem}\|_{H^\alpha}^2+\l M_0(2\mu+\nu)\|\u_{\rm
rem}\|_{H^{\alpha+1}}^2.
\end{align}

We deal with the fourth term on the right-hand side of
\eqref{higher1}. From \eqref{rem2}, we can easily get the equation
for $D^\alpha\rho_{\rm rem}$,
\begin{align}
\partial_t D^\alpha\rho_{\rm rem}+\u_\lambda\cdot\nabla
D^\alpha\rho_{\rm
rem}+\frac{1}{\lambda}\rho_\lambda\mbox{div}D^\alpha\u_{\rm
rem}=D^\alpha h_0+h^\alpha \label{rhoa}
\end{align}
with
\begin{align}
h^\alpha=&-D^\alpha(\u_\lambda\cdot\nabla\rho_{\rm rem})+\u_{\l
 }\cdot\nabla D^\alpha
\rho_{\rm rem}-\frac{1}{\lambda}D^\alpha(\rho_\lambda\mbox{div}\u_{\rm rem})\nonumber\\
&+\frac{1}{\lambda}\rho_\lambda\mbox{div}D^\alpha\u_{\rm
rem}.\nonumber
\end{align}
In view of \eqref{rhoa} and the Poisson equation \eqref{rem4}$_2$,
we get
\begin{align}
&\dfrac{2}{\l}\dint \langle A_0(U_{\rm rem}) D^\alpha B,D^\alpha
U_{\rm rem} \rangle dx=-\dfrac{2}{\l}\dint \rho_\lambda\nabla D^\alpha\phi_{\rm rem} D^\alpha\u_{\rm rem}dx\nonumber\\
=&\dfrac{2}{\l}\dint \rho_\lambda\mbox{div}D^\alpha\u_{\rm
rem}D^\alpha\phi_{\rm rem}dx+\dfrac{2}{\l}\dint
\nabla\rho_\lambda D^\alpha\u_{\rm rem}D^\alpha\phi_{\rm rem}dx\nonumber\\
=&-2\dint\partial_t D^\alpha\rho_{\rm rem}D^\alpha\phi_{\rm rem}
dx-2\dint\u_\lambda\nabla D^\alpha\rho_{\rm rem}D^\alpha\phi_{\rm
rem} dx \nonumber\\&+2\dint D^\alpha h_0 D^\alpha\phi_{\rm rem}
dx+2\dint \nabla(\rho_{\rm osc}+\lambda(\Delta
\Pi+\rho_{\rm cor})+\lambda\rho_{\rm rem})D^\alpha\u_{\rm
rem}D^\alpha\phi_{\rm rem} dx \nonumber\\& +2\dint h^\alpha
D^\alpha\phi_{\rm rem} dx
\nonumber\\
\leq &-\dfrac{d}{dt}\|D^\alpha \nabla\phi_{\rm rem}\|_{L^2}^2
+M_0(1+\lambda M)\bigg(\|D^\alpha\nabla\phi_{\rm
rem}\|_{L^2}^2+
\dsum_{0\leq|\beta|\leq|\alpha|}\|D^{\beta}U_{\rm rem}\|_{E}^2\bigg) \nonumber\\
 &\quad +\epsilon_3 \dint \|\nabla D^{\alpha}\u_{\rm rem}\|^2 dx
\end{align}
for some sufficiently small constant $\epsilon_3>0$.

The fifth term on the right-hand side of \eqref{higher1} is very
tedious. The main techniques involved are   Leibniz's formula,
Moser-type calculus inequalities \eqref{MI1}-\eqref{MI2}, and
Sobolev's embedding inequalities.  Actually, after the tedious
computations, we finally obtain the following estimate
\begin{align}
&2\dint\langle A_0(U_{\rm rem})D^\alpha F(x,t,  U_{\rm rem}),
D^\alpha
U_{\rm rem}\rangle dx\nonumber\\
\leq&
\lambda (2\mu+\nu+\kappa)M_0(1+M)\bigg[ \dsum_{0\leq|\beta|\leq
|\alpha|}\big(\|\nabla D^\beta \u_{\rm rem}\|_{L^2}^2
+\|\nabla D^\beta \theta_{\rm rem}\|_{L^2}^2\big)\nonumber\\
&\qquad\qquad +\dsum_{0\leq|\beta|\leq |\alpha|}
\|D^\beta U_{\rm rem}\|_{E}^2 \bigg]
+(2\mu+\nu+\kappa+1)M_0.\label{F2}
\end{align}

The commutating term $H^\alpha $ can be bounded by
\begin{align}
&\dint\langle H^\alpha(t), D^\alpha U_{\rm rem} \rangle dx\nonumber\\
\leq & \dsum_{1\leq |\beta|\leq |\alpha|}M_0(1+\lambda M)\|D^\beta
U_{\rm rem}\|_{E}^2+\|D^\alpha U_{\rm rem}\|_{E}^2+M_0.\label{H1}
\end{align}

We now re-choose $\delta_T$ sufficiently  small such that, for
$\lambda\in (0, \delta_T]$,
\begin{align}
\lambda sM_0(M+1)(2\mu+\nu+\kappa)\leq \min\Big\{\frac{\xi_2}{2},
\frac{\kappa\rho_-}{2\theta_+}\Big\}:=\eta_2. \label{lamsmall}
\end{align}
Let
\begin{align}
\Phi(t)=\lambda^2\|\rho_{\rm rem}\|_{H^s}^2+\|\u_{\rm
rem}\|_{H^s}^2+\|\theta_{\rm rem}\|_{H^s}^2.
\end{align}
Taking   $\delta$  and $\epsilon_3$  small enough and combining the estimates
\eqref{Gamma}-\eqref{H1} with \eqref{higher1} and
\eqref{estimatezero}, we obtain that
\begin{align}
&c_0\Phi(t)+\|\nabla\phi_{\rm
rem}\|_{H^s}^2+\frac{\xi}{2}\dint_0^t\|\u_{\rm rem}\|_{H^{s+1}}^2
dr+ \frac{\kappa\rho_-}{2\theta_+}\dint_0^t\|\theta\|_{H^{s+1}}^2dr\nonumber\\
\leq &
\dint_0^t\Big\{M_0\Big(M_0(1+\l(M+\tilde{M}))+3\eta+M_0\kappa+\lambda(2\mu+\nu)
MM_0\|\u_{\rm rem}\|_{H^{s+1}}^{\frac{3}{2}}\Big)\nonumber\\&
\times\Big(c_0\Phi(r)+\|\nabla\phi_{\rm
rem}\|_{H^s}^2(r)\Big)\Big\} dr+ c_0\Phi(0)+\|\nabla\phi_{\rm
rem}(0)\|_{H^s}^2+M_0(2\mu+\nu+\kappa)T,\label{integral}
\end{align}
where $\xi=\min\{\xi_1,\xi_2\}$ and $\eta=\max\{\eta_1,\eta_2\}$.
By virtue of Gronwall's inequality, we obtain that
\begin{align}
&c_0\Phi(t)+\|\nabla\phi_{\rm rem}\|_{H^s}^2 \leq
\big(c_0\Phi(0)+\|\nabla\phi_{\rm
rem}(0)\|_{H^s}^2+M_0(2\mu+\nu+\kappa)T\big)\nonumber\\
&\times{\rm exp}\Big\{
M_0\dint_0^t \Big[ M_0(1+\l(M+\tilde{M}))+3\eta+M_0\kappa+\lambda(2\mu+\nu)
MM_0\|\u_{\rm rem}\|_{H^{s+1}}^{\frac{3}{2}}\Big]dr\Big\}.\label{Gr1}
\end{align}
From \eqref{apriori1} and H\"older's inequality, we have
\begin{align} \lambda(2\mu+\nu) MM_0\dint_0^t\|\u_{\rm
rem}\|_{H^{s+1}}^{\frac{3}{2}}dr\leq \lambda
M_0(2\mu+\nu)M^{\frac{7}{4}}T^{\frac{1}{4}}.\label{23}
\end{align}
In view of \eqref{assin} and \eqref{assin1}, we obtain that
\begin{align}
\lambda^2\|\rho_{\rm rem}(0)\|_{H^s}^2\leq \tilde{C}\lambda^2,\,\,
\|\u_{\rm rem}(0)\|_{H^s}^2+\|\theta_{\rm rem}(0)\|_{H^s}^2\leq
\tilde{C}\label{ie1}
\end{align}
and
\begin{align}
\|\nabla\phi_{\rm rem}\|_{H^s}^2\leq \tilde{C}.\label{ie2}
\end{align}
We choose  $\delta_T$ sufficiently small such that,   for $\lambda
\in (0, \delta_T]$, it holds that
\begin{align}
\l(M+\tilde{M})+\l (2\mu+\nu)M^{\frac{7}{4}}<1.\label{small}
\end{align}
Set
\begin{align*}
L_1=M_0(2M_0+3\eta+M_0\kappa+M_0T^{1/4}).
\end{align*}
Substituting \eqref{23}-\eqref{small} into \eqref{Gr1}, we obtain
that
\begin{align}
c_0\Phi(t)+\|\nabla\phi_{\rm rem}\|_{H^s}^2 \leq&
(c_0\Phi(0)+\|\nabla\phi_{\rm rem}(0)\|_{H^s}^2+M_0(2\mu+\nu+\kappa)T)e^{L_1T}\nonumber\\
\leq& (M_0\tilde{C}+M_0(2\mu+\nu+\kappa)T)e^{L_1T}=:L_3.
\end{align}
In view of \eqref{integral}, we get that
\begin{align}
\frac{\xi}{2}\dint_0^t\|\u_{\rm rem}\|_{H^{s+1}}^2 dr+
\frac{\kappa\rho_-}{2\theta_+}\dint_0^t\|\theta\|_{H^{s+1}}^2dr\leq
L_1 L_3T +M_0\tilde{C}+M_0(2\mu+\nu+\kappa)T.
\end{align}
Therefore \eqref{apriori1} is proved if we set
\begin{align}
M^2=:(L_3+L_1 L_3T
+M_0\tilde{C}+M_0(2\mu+\nu+\kappa)T)\cdot\max\Big\{\frac{1}{c_0},
1, \frac{2}{\xi}, \frac{2\theta_+}{\kappa\rho_-}\Big\}.\label{M1}
\end{align}

It follows from \eqref{rem4}  that \begin{align} &\dsup_{0\leq t\leq
T}\Big(\l\|\partial_t\rho_{\rm rem}(t)\|_{H^{s-1}}+
\l\|\partial_t\u_{\rm rem}(t)\|_{H^{s-2}}+
\|\partial_t\theta_{\rm rem}(t)\|_{H^{s-1}} \nonumber\\
&\qquad\qquad +\l\|\partial_t\nabla\phi_{\rm
rem}(t)\|_{H^s}\Big)\leq \tilde{M}
\end{align}
with
\begin{align}
\tilde{M}:=(M_0(1+2M))^{1/2}.\label{M2}
\end{align}
The proof of Lemma \ref{apriori} is completed.

\end{proof}

\begin{proof}[Proof of Theorem \ref{Thm3}]  With the a priori estimates \eqref{apriori1} and \eqref{apriorit},
we now start the proof of Theorem \ref{Thm3}. We first construct
the approximate solutions. Define
\begin{align*}
(U_{\rm rem}^{n+1}, \phi_{\rm rem}^{n+1})=(\rho_{\rm rem}^{n+1},
\u_{\rm rem}^{n+1}, \theta_{\rm rem}^{n+1}, \phi_{\rm
rem}^{n+1})^{\rm{T}}\quad  (  n\geq 0)
\end{align*}
inductively as the solution of linear equations
\begin{equation}\label{appro1}
\left\{
\begin{aligned} &A_0(U_{\rm rem}^n)\partial_tU_{\rm
rem}^{n+1}+\dsum_{j=1}^N\mathcal{A}_j(x,t,U_{\rm rem}^n)
\partial_{x_j}U_{\rm rem}^{n+1}-
\mu\rho_\l^n\Delta\tilde{\u}_{\rm rem}^{n+1} \\
& \qquad-(\mu+\nu)\rho_\l^n\nabla\mbox{div}\tilde{\u}_{\rm
rem}^{n+1}
-\frac{\kappa\rho_\lambda^n}{\theta_\lambda^n}\Delta \tilde{\theta}_{\rm rem}^{n+1}  \\
&\qquad =\lambda\nu \tilde{J}^{n}+\frac{\l\mu}{2}
\tilde{G}^n+\dfrac{1}{\l}\tilde{B}^{n+1}+
 \tilde{F}^{n}, \\
&-\Delta \phi^n_{\rm rem}=\rho^n_{\rm rem}
\end{aligned}\right.
\end{equation}
with initial data
\begin{equation}\label{ini5}
  U^n_{\rm rem}(x,0)=U_{{\rm rem} 0}(x),
\end{equation}
 where
\begin{align}
&\rho_\l^n(x,t)=1+\l\rho_{\rm osc}(x,t)+\l^2(\Delta
\Pi(x,t)+\rho_{\rm cor}(x,t/\l))+\l^2\rho_{\rm rem}^n(x,t), \nonumber\\
&\u^n_{\l}(x,t)=\v+\u_{\rm osc}(x,t)+\l\u_{\rm cor}(x,t/\l)+\l\u^n_{\rm rem}(x,t),\nonumber\\
&\theta^n_{\l}(x,t)=\theta(x,t)+\l \theta_{\rm cor}(x,t/\l)+\l
\theta^n_{\rm rem}(x,t),\nonumber\\
&\phi^n_\l(x,t)=\phi_{\rm osc}(x,t)+\l(\Pi(x,t)+\phi_{\rm
cor}(x,t/\l))+\l\phi^n_{\rm rem}(x,t),\nonumber\\
&\tilde{\u}_{\rm rem}^{n+1}=(0, \u_{\rm rem}^{n+1},
0)^{\rm{T}},\quad \tilde{B}^{n+1}=A_0B(x,t, U_{\rm rem})=
(0, -\rho_\l^n\nabla\phi_{\rm rem}^{n+1}, 0),\nonumber\\
&\tilde{J}^n:=A_0D(x,t, U_{\rm
rem}^n)=\Big(0,\dots,0,\frac{\rho_\l^n }{
\theta_\l^n}(\mbox{div}\u_{\rm rem}^n)^2\Big)^{\rm{T}},\nonumber\\
&\tilde{G}^n:=A_0G(x,t,  U_{\rm
rem}^n)=\Big(0,\dots,0,\frac{\rho_\l^n }{
\theta_\l^n}\sum^N_{i,j=1}((\partial_i u_{\rm rem}^{j})^n+(\partial_ju_{\rm rem}^{i})^n)^2\Big)^{\rm{T}},\nonumber\\
&\tilde{F}^{n}=A_0F(x,t, U_{\rm rem}^{n}).\nonumber
\end{align}

It is standard to know that the approximate problem \eqref{appro1}
admits a unique solution such that \begin{align*} &(\rho_{\rm
rem}^{n+1}, \u_{\rm rem}^{n+1}, \theta_{\rm rem}^{n+1},
\nabla\phi_{\rm rem}^{n+1})\in C([0,T];H^s),&& \nabla\phi_{\rm
rem}^{n+1}\in C([0,T];H^{s+1}), \\
&   \u_{\rm rem}^{n+1}\in L^2(0,T;H^{s+1}), &&  \theta_{\rm
rem}^{n+1}\in L^2(0,T;H^{s+1}), \\
 &  \partial_t\rho_{\rm
rem}^{n+1}\in C([0,T];H^{s-1}), && \partial_t\u_{\rm rem}^{n+1}\in C([0,T];H^{s-2}), \\
&    \partial_t \theta_{\rm rem}^{n+1}\in C([0,T];H^{s-2}), &&
\partial_t\nabla\phi_{\rm rem}^{n+1}\in C([0,T];H^s),
\end{align*}
and satisfies the uniform estimates
\begin{align}
&\dsup_{0\leq t\leq T}(\|(\lambda\rho_{\rm rem}^{n+1}, \u_{\rm
rem}^{n+1},
\theta_{\rm rem}^{n+1})(t)\|_{H^s}^2+\|\nabla\phi_{\rm rem}^{n+1}\|_{H^{s+1}}^2),\nonumber\\
&\qquad+\dint_0^T \|\u_{\rm rem}^{n+1}\|_{H^{s+1}}^2 dt
+\dint_0^T \|\theta_{\rm rem}^{n+1}\|_{H^{s+1}}^2 \;dt \leq M^2, \label{apriori1r}\\
&\dsup_{0\leq t\leq T}\Big(\l^2\|\partial_t\rho_{\rm
rem}^{n+1}(t)\|_{H^{s-1}}^2+ \l^2\|\partial_t\u_{\rm
rem}^{n+1}(t)\|_{H^{s-2}}^2 +\|\partial_t \theta_{\rm
rem}^{n+1}(t)\|_{H^{s-1}}^2\nonumber\\
&\qquad+\l^2\|\partial_t\nabla\phi_{\rm
rem}^{n+1}(t)\|_{H^s}^2\Big)\leq \tilde{M}^2.
\end{align}
 It is standard to verify that the difference
\begin{align*} (\bar{\rho}_{\rm rem}^{n+1},
\bar{\u}_{\rm rem}^{n+1}, \bar{\theta}_{\rm rem}^{n+1},
\bar{\phi}_{\rm rem}^{n+1})= (\rho_{\rm rem}^{n+1}-\rho_{\rm
rem}^n, \u_{\rm rem}^{n+1}-\u_{\rm rem}^n, \theta_{\rm
rem}^{n+1}-\theta_{\rm rem}^n, \phi_{\rm rem}^{n+1}-\phi_{\rm
rem}^n)
\end{align*}
satisfies
\begin{equation}\label{diff}
\left\{\begin{aligned}
 &\partial_t\bar{\rho}_{\rm
rem}^{n+1}+\u_\l^n\nabla\bar{\rho}_{\rm rem}^{n+1}
+\frac{1}{\l}\rho_\l^n\mbox{div}\bar{\u}_{\rm rem}^{n+1} \\
&\quad =-\l\bar{\u}_{\rm rem}^{n}\nabla\rho_{\rm rem}^n-
\l\bar{\rho}_{\rm rem}^{n}\mbox{div}\u_{\rm rem}^n \\
&\qquad +h_0(x,t,\u_{\rm rem}^{n+1}, \rho_{\rm rem}^{n+1})-h_0(x,t,\u_{\rm
rem}^{n},
\rho_{\rm rem}^{n}), \\
&\partial_t\bar{\u}_{\rm
rem}^{n+1}+(\u_\l^n\cdot\nabla)\bar{\u}_{\rm rem}^{n+1} +\l
\frac{\theta_\l^n}{\rho_\l^n}\nabla\bar{\rho}_{\rm
rem}^{n+1}+\nabla\bar{\theta}_{\rm rem}^{n+1} \\
&\quad -\mu\Delta\bar{\u}_{\rm rem}^{n+1}
-(\mu+\nu)\nabla\mbox{div}\bar{\u}_{\rm rem}^{n+1} \\
&\qquad =\frac{1}{\l}\nabla\phi_{\rm rem}^{n+1}-\lambda(\bar{\u}_{\rm
rem}^n\cdot\nabla)\u_{\rm rem}^n
-\l\Big(\frac{\theta_\l^n}{\rho_\l^n}-\frac{\theta_\l^{n-1}}{\rho_\l^{n-1}}\Big)\nabla\rho_{\rm rem}^n \\
&\quad \qquad+f_0(x,t,\u_{\rm rem}^{n+1}, \rho_{\rm rem}^{n+1})
-f_0(x,t,\u_{\rm rem}^{n}, \rho_{\rm rem}^{n}), \\
&\partial_t\bar{\theta}_{\rm
rem}^{n+1}+\u_\l^n\cdot\nabla\bar{\theta}_{\rm rem}^{n+1}+
\theta_\l^n\mbox{div}\bar{\u}_{\rm rem}^{n+1} \\
&\quad =\l\nu(\tilde{J}^n-\tilde{J}^{n-1})+2\l\mu(\tilde{G}^n-\tilde{G}^{n-1})-(
\theta_\l^n-
\theta_\l^{n-1})\mbox{div}\u_{\rm rem}^n \\
&\qquad -\l\bar{\u}_{\rm rem}^n \nabla \theta_{\rm rem}^n+g_0(x,t,\u_{\rm
rem}^{n+1}, \rho_{\rm rem}^{n+1})-g_0(x,t,\u_{\rm rem}^{n},
\rho_{\rm rem}^{n}).
\end{aligned}\right.
\end{equation}
Observing that, for $|\alpha|\leq s$,
\begin{align}
&|D^\alpha(\tilde{J}^n-\tilde{J}^{n-1})|+|D^\alpha(\tilde{G}^n-\tilde{G}^{n-1})|\nonumber\\
\leq & M_0\dsum_{|\alpha|-1=|\beta|+|\gamma|\leq s-1}\big[(|D^{\beta+1}\u_{\rm
rem}^n|+|D^{\beta+1}\u_{\rm rem}^n|)|D^{\gamma+1}\bar{\u}_{\rm
rem}^n|\big].
\end{align}

Then repeating the previous analysis used in the proof of Lemma
\ref{apriori} and using the interpolation inequalities, we can show
that there is a $\delta_T>0$ such that, for any $\l\in (0,
\delta_T]$ and $s'<s$,
\begin{align}
&\dsup_{0\leq t\leq T}\big(\|(\lambda\bar{\rho}_{\rm
rem}^{n+1},\bar{\u}_{\rm rem}^{n+1},\bar{\theta}_{\rm
rem}^{n+1})(t)\|_{H^{s'}}^2 +\|\nabla\bar{\phi}_{\rm
rem}^{n+1}(t)\|_{H^{s'+1}}^2\big)\nonumber\\
&\qquad +\dint_0^T\|\bar{\u}_{\rm rem}^{n+1}\|_{H^{s'+1}}^2
dr+\dint_0^T\|\bar{\theta}_{\rm
rem}^{n+1}\|_{H^{s'+1}}^2  \leq C,\nonumber\\
&\dsup_{0\leq t\leq T}\big(\lambda^2\|\partial_t\bar{\rho}_{\rm
rem}^{n+1}(t)\|_{H^{s'-1}}^2+\lambda^2\|\partial_t\bar{\u}_{\rm
rem}^{n+1}(t)\|_{H^{s'-2}}^2+\|\partial_t\bar{\theta}_{\rm
rem}^{n+1}
(t)\|_{H^{s'-2}}^2\nonumber\\
&\qquad+\lambda^2\|\partial_t\nabla\bar{\phi}_{\rm
rem}^{n+1}(t)\|_{H^{s'}}^2\big) \leq C\nonumber
\end{align}
for some constant $C>0$. Then the Arzel\`{a}-Ascoli theorem
  implies that there exists a limit vector
function
\begin{align}
(\rho_{\rm rem}, \u_{\rm rem}, \theta_{\rm rem}, \nabla\phi_{\rm
rem})^{\rm T}\in L^\infty(0,T; H^{s'})\cap{\rm Lip}([0,T];
H^{s'-1})\nonumber
\end{align}
satisfying \eqref{apriori1}-\eqref{apriorit} such  that
\begin{align}
\sup_{0\leq t\leq T}\|(\rho_{\rm rem}^{n+1}-\rho_{\rm rem},
\u_{\rm rem}^{n+1}-\u_{\rm rem}, \theta_{\rm
rem}^{n+1}-\theta_{\rm rem},\nabla\phi_{\rm
rem}^{n+1}-\nabla\phi_{\rm rem})(t)\|_{H^{s'-2}}\rightarrow
0\nonumber
\end{align}
as $n\rightarrow +\infty$  for any $\lambda\in (0, \delta_T]$.
Furthermore, for $N/2-[N/2]<\sigma<1$, we have the convergence
\begin{align}
(\rho_{\rm rem}^{n+1},\u_{\rm rem}^{n+1},\theta_{\rm rem}^{n+1},
\nabla\phi_{\rm rem}^{n+1})^{\rm T}\rightarrow(\rho_{\rm rem},
\u_{\rm rem}, \theta_{\rm rem}, \nabla\phi_{\rm rem})^{\rm T}\nonumber
\end{align}
in $C([0, T]; H^{s-\sigma})$ by the standard interpolation
inequality. Moreover, by Sobolev's embedding theorem, we have
\begin{align}
&(\rho_{\rm rem}, \u_{\rm rem}, \theta_{\rm rem}, \phi_{\rm
rem})^{\rm T}\in C([0,T]; H^{s'})\cap C^1([0,T];
H^{s'-2})\nonumber\\
&\qquad\qquad \hookrightarrow C^1([0,T]\times \mathbb{T}^N)\cap
C([0,T]; C^2(\mathbb{T}^N))\nonumber
\end{align}
for any $\lambda\in (0, \delta_T]$, where we have used the fact
  $s'>N/2+2$.  Then the existence of classical solutions to the
initial value problem \eqref{rem2}-\eqref{ini2} is proved. The
uniqueness of the classical solutions can be proved easily by
energy estimates for the difference of any two solutions. Thus the
proof of Theorem \ref{Thm3} is finished.

\end{proof}

\section{Proofs of Theorem \ref{Thm} and  Theorem \ref{Thm2}}

\begin{proof}[Proof of Theorem \ref{Thm}]
By the asymptotic expansion \eqref{expa}, Propositions \ref{Proa}
and \ref{Prob}, the existence and uniqueness of classical solutions
to the initial value problem of Navier-Stokes-Poisson system
\eqref{pns4}-\eqref{pns7} is proved and the solution satisfies
\begin{align}
&\sup_{0\leq t\leq T}\|(\rho_\lambda, \u_\lambda, \theta_\lambda)
(t)\|_{H^s}+\sup_{0\leq t\leq T}\|\nabla\phi_\lambda (t)\|_{H^{s+1}}\nonumber\\
&\qquad\qquad +\|\u_\lambda\|_{L^2(0,T;
H^{s+1})}+\|\theta_\lambda\|_{L^2(0,T; H^{s+1})}\leq
C(T), \nonumber\\
&\dsup_{0\leq t\leq T}\big(\|\partial_t( \rho_\lambda, \u_\lambda,
\theta_\lambda)(t)\|_{H^s}+\|\partial_t\nabla\phi_\lambda(t)\|_{H^{s+1}}\big)
\leq C(T, \lambda), \nonumber
\end{align}
where $C(T)>0$ is a constant independent of $\lambda$ and $C(T,
\lambda)>0$ is a constant dependent on $\lambda$. Moreover, it is
easy to see that, for $\lambda\in (0, \delta_T]$,
\begin{align}
& \sup_{0\leq t\leq T}\|(\rho_\lambda-1, \u_\lambda-\v-\u_{\rm osc},
\theta_\lambda-\theta)(t)\|_{H^s}\nonumber\\
&\qquad +\sup_{0\leq t\leq T}\|(\nabla\phi_\lambda-\nabla\phi_{\rm
osc} )(t)\|_{H^{s+1}}\leq C(T)\lambda\nonumber.
\end{align}
Thus the proof of Theorem \ref{Thm} is finished.
\end{proof}

As far as the combined quasineutral, vanishing viscosity and
vanishing heat conductivity limit is concerned, we can follow the
same lines as the proof of Theorem \ref{Thm}. Recalling the
uniformly bounded estimates obtained in Lemma \ref{apriori}, we are
able to get the uniform bound with respect to $\lambda, \mu, \nu $
and $\kappa$ for the solutions. Thus Theorem \ref{Thm2} can be
proved similarly with minor modifications of our previous arguments.
We omit the details here for conciseness.

\medskip
{\bf Acknowledgements}\ \  Ju is supported   by   NSFC (Grant
10701011).   F. Li is  supported   by NSFC (Grant   10501047).  H.
Li is supported  by NSFC (Grants 10431060, 10871134), the Beijing
Nova program, the NCET support of the Ministry of Education of
China, the Huo Ying Dong Foundation 111033, the support of Institute
of Mathematics and Interdisciplinary Science at CNU.

\end{document}